\documentclass[a4paper]{amsart}
\usepackage{amssymb}

\newcommand{\heute}{27 October 1999}

\newtheorem{theorem}{Theorem}[section]
\newtheorem{lemma}[theorem]{Lemma}
\newtheorem{corollary}[theorem]{Corollary}
\newtheorem{proposition}[theorem]{Proposition}

\theoremstyle{definition}

\theoremstyle{remark}
\newtheorem{remark}[theorem]{Remark}

\newcommand{\z}{\mathbb{Z}}

\newcommand{\f}[1][p]{\mathbb{F}_{#1}}
\newcommand{\Aut}{\operatorname{Aut}}
\newcommand{\Bild}{\operatorname{Im}}

\newcommand{\Out}{\operatorname{Out}}
\newcommand{\Res}{\operatorname{Res}}
\newcommand{\Cor}{\operatorname{Cor}}

\newcommand{\GL}{\text{\sl GL}}
\newcommand{\GLn}[1][p]{\GL(\f)}

\begin{document}

\title{Almost all extraspecial $p$-groups are {S}wan groups}

\author[D.~J. Green]{David John Green}
\address{Dept of Mathematics \\ Univ.\@ of Wuppertal \\
D--42097 Wuppertal \\ Germany}
\email{green@math.uni-wuppertal.de}

\author[P.~A. Minh]{Pham Anh Minh}
\address{Dept of Mathematics \\ University of Hue \\ 27 Nguyen Hue \\
Hue \\ Vietnam}

\date{\heute}

\begin{abstract}
Let $P$ be an extraspecial $p$-group which is neither dihedral of order~$8$,
nor of odd order~$p^3$ and exponent~$p$. Let $G$ be a finite group having
$P$ as a Sylow $p$-subgroup. Then the mod-$p$ cohomology ring
of~$G$ coincides with that of the normalizer $N_G(P)$.
\end{abstract}

\maketitle

\section*{Introduction}
Let~$P$ be a finite $p$-group.  Martino and Priddy call $P$~a Swan
group~\cite{MartinoPriddy:Swan} if for every finite group~$G$ with Sylow
$p$-subgroup~$P$, the mod-$p$ cohomology ring $H^*(G)$ coincides with
$H^*(N_G(P)) \cong H^*(P)^{N_G(P)}$.
In particular, if there are no so-called transfer summands in the stable
decomposition of the classifying space~$BP$, then $P$~is a Swan group.

We prove in Theorem~\ref{theorem:main} that all extraspecial $p$-groups
are Swan groups, apart from the well-known exceptions $2^{1+2}_+ = D_8$ and
(for $p$~odd) $p^{1+2}_+ = E$.
The cases where $P$~is the metacyclic group $p^{1+2}_- = M(p^3)$ with $p$~odd,
and where $P$~is $2^{1+2n}_- = Q_8 * D_8 * \cdots * D_8$,
were proved in~\cite{MartinoPriddy:Swan};
the former case being due to G.~Glauberman.
Earlier, the $p=2$ case of the theorem was published
in~\cite{Ogawa:extraspecial}, but with an incorrect proof:
see Remark~\ref{remark:Ogawa}.
In Corollary~\ref{coroll:transfer} we generalize another result of Martino
and Priddy, exhibiting three infinite families of Swan groups whose
classifying spaces do have transfer summands in their stable decompositions.

Throughout this paper we denote the mod-$p$ cohomology ring $H^*(G, \f[p])$
by $H^*(G)$.  A suitable reference on group cohomology
is Evens' book~\cite{Evens:book}.

\section{Extraspecial $p$-groups}
Recall that a $p$-group~$P$ is called extraspecial if its centre $Z(P)$,
its derived subgroup~$P'$ and its Frattini subgroup~$\Phi(P)$ all coincide,
and are cyclic of order~$p$.  So if $P$~is extraspecial there is
a central extension
\[
1 \longrightarrow \f \longrightarrow P \stackrel{\psi}{\longrightarrow} V
\longrightarrow 1
\]
with $V$~an elementary abelian $p$-group.
Hence there is a nondegenerate alternate bilinear form $f$~on $V$
defined by
$f(\psi(g), \psi(h)) = [g,h]$ for all $g,h \in P$.

Moreover, for $p$~odd there is a linear form $\lambda$~on $V$ defined
by $\lambda(\psi g) = g^p$; and for $p=2$ there is a quadratic form
$Q$~on $V$ defined by $Q(\psi g) = g^2$ with associated bilinear form~$f$.
Conversely, such a pair $(\lambda, f)$ determines an extraspecial $p$-group
when $p$~is odd, and such a~$Q$ determines an extraspecial $2$-group.

See \cite[Ch.~6]{Jacobson:BasicI} for a reference on alternate
bilinear forms, and \cite[I.16]{Dieudonne:3rdEd} for a reference on quadratic
forms in characteristic~$2$.  Nondegeneracy means that $V$~has even
dimension, say~$2n$.  Up to change of basis for~$V$ there are two
possibilities for $Q$~when $p=2$, and exactly one choice of~$f$ for odd~$p$.
For $p$~odd there are two possibilities for
the pair $(\lambda,f)$: either $\lambda$~is identically
zero, or it is not.  Note that in the case of nonzero~$\lambda$, Witt's
extension theorem does not hold for the pair $(\lambda, f)$, as the
restriction of $f$~to $\ker(\lambda)$ does have a kernel.

So there are four types of extraspecial $p$-groups. In each case
$P$ is generated by $A_1,\ldots,A_n,B_1, \ldots,B_n, C$, with $C$~central
of order~$p$, $[A_i, A_j] = [B_i, B_j] = 1$ and $[A_i, B_j] = C^{\delta_{ij}}$.
Moreover $A_i^p = B_i^p = 1$ for $2 \leq i \leq n$.
The four cases are:
\begin{itemize}
\item $2^{1+2n}_+ = D_8 * \cdots * D_8$:
here $A_1^2 = B_1^2 = 1$.
\item $2^{1+2n}_- = Q_8 * D_8 * \cdots * D_8$:
here $A_1^2 = B_1^2 = C$.
\item $p^{1+2n}_+ = E * \cdots * E$ has odd exponent~$p$:
here $A_1^p = B_1^p = 1$.
\item $p^{1+2n}_- = M(p^3) * E * \cdots * E$ has odd exponent~$p^2$:
here $A_1^p = C$ and $B_1^p = 1$.
\end{itemize}

\noindent
The characteristic subgroup $\Omega_1(P)$~of $P$ is the
subgroup generated by all order~$p$ elements.
Denote $Z(\Omega_1(P))$ by~$Y$.  If $P$~has odd exponent~$p^2$
then $Y = \langle B_1, C \rangle$ is rank two elementary abelian; in all other
cases, $Y$~equals $Z = Z(P)$.

The following result could be called Witt's theorem for
extraspecial $p$-groups.

\begin{proposition}
\label{prop:Witt}
Let $P$ be an extraspecial $p$-group. Suppose that $H,K$ are subgroups
of~$P$ containing~$Z$, and that $\phi \colon H \rightarrow K$ is a
group isomorphism inducing the identity map on~$Z$.
If $P$~has odd exponent~$p^2$, assume further that
$H \cap Y = K \cap Y$ and that $\phi$ induces the identity map
on $(H \cap Y) / Z$.
Then $\phi$~extends to an automorphism of~$P$.
\end{proposition}

\begin{lemma}
\label{lemma:Wittprep}
In Proposition~\ref{prop:Witt}, suppose that $P$~has odd exponent~$p^2$
and that $H \cap Y$ is $Z$ rather than~$Y$. Then $\phi$~extends to an
isomorphism from $\langle H, Y \rangle$ to $\langle K, Y \rangle$ which
itself satisfies the conditions of Proposition~\ref{prop:Witt}.
\end{lemma}

\begin{proof}
Since $\phi(C) = C$ it follows that $h^{-1} \phi(h)$ lies in $\ker(\lambda)$
for every $h \in H$.
Hence $\phi([h, B_1]) = [\phi(h), B_1]$. So we may set $\phi(B_1) = B_1$.
\end{proof}

\begin{proof}[Proof of Proposition~\ref{prop:Witt}]
Denote by $U,W$ the images in~$V$ of $H,K$ respectively.  Since $\phi$~is
the identity on~$Z$ there is an $\f$-vector space isomorphism
$\rho \colon U \rightarrow W$ induced by~$\phi$ which respects the
alternate bilinear form $f$~on $V$.

If $p$~is $2$ then $\rho$ respects the quadratic form~$Q$.  Since
Witt's extension theorem holds for~$Q$ (see~\cite[p.~36]{Dieudonne:3rdEd}),
we may extend $\rho$ to a $Q$-orthogonal transformation $\rho$~of $V$.
Using the standard generators for~$P$ we may lift $\rho$~to an
automorphism $\phi'$ of~$P$. If $h \in H$ then $\phi'(h) = \phi(h) C^r$ for
some $r \in \z/p$.  Since $P$ has enough inner automorphisms,
we may assume that $\phi'$~extends $\phi$.

To be more precise: pick $h_1,\ldots,h_m \in H$ whose images under
$\psi \colon P \rightarrow V$ constitute a basis for~$U$.
Since the alternate bilinear form $f$~on $V$ is nondegenerate we can pick
$g_1,\ldots,g_m \in P$ such that $f(\psi(g_i), \psi(\phi h_j)) = \delta_{ij}$.
Hence conjugation by~$g_i$ fixes $\phi(h_j)$ for $j \not = i$ and
sends $\phi(h_i)$ to $\phi(h_i) C$.  So we can correct $\phi'$ by an inner
automorphism of~$P$ to ensure that $\phi' = \phi$ on $H$.

Now suppose that $p$~is odd.  As Witt's extension theorem holds for~$f$
(see \cite[6.9]{Jacobson:BasicI}), we may extend~$\rho$ to a transformation
of~$V$ that respects~$f$.
But then $\rho$~respects $\lambda$ too: this is trivial in the exponent~$p$
case, as~$\lambda$ is then zero. In the exponent~$p^2$ case,
we may assume by Lemma~\ref{lemma:Wittprep} that $H,K$ contain~$Y$
and that $\rho$ fixes $\psi(B_1)$. So $\rho$ respects $\lambda$ by
Lemma~\ref{lemma:lambda} below.
As in the $p=2$ case we can now lift $\rho$ to an automorphism of~$P$ which
extends~$\phi$.
\end{proof}

\begin{lemma}
\label{lemma:lambda}
Suppose $P$ is an extraspecial $p$-group with odd exponent~$p^2$.
Taking $C$ as basis for $Z \cong \f$, we may assume that
$\lambda(\psi A_1) = 1$.  Then for any $v \in V$ we have
$\lambda(v) = f(v, \psi(B_1))$.
\end{lemma}

\begin{proof}
Each element $g$~of $P$ has canonical form $B_1^{s_1} \cdots B_n^{s_n} \cdot
A_1^{r_1} \cdots A_n^{r_n} \cdot C^t$.  Then
$g^p = C^{r_1} = [g, B_1]$.
\end{proof}

\section{Local subgroup structure}
Throughout this section $G$ is a finite group with extraspecial
Sylow $p$-subgroup~$P$.

\begin{lemma}
\label{lemma:centralizerFrattini}
Suppose that $P$~is not one of $D_8$, $E$, $M(p^3)$.
Then for any order~$p$ element $g$~of $P$,
the centralizer of $g$~in $P$ has the same Frattini subgroup as $P$~itself.
\end{lemma}

\begin{proof}
If $n \geq 2$ then $C_P(g)$ is not abelian. In $Q_8$, all order~$p$ elements
are central. The other three groups really are exceptions: take $g = B_1$.
\end{proof}

\begin{lemma}
\label{lemma:Z}
Suppose that the centralizer of each exponent~$p$ element of~$P$ has~$Z$ as
its Frattini subgroup.
Then $Z^g = Z$ for every $g \in G$ such that $P \cap P^g$ contains~$Z$.
Moreover, such $g$~factorize as $g = g_1 g_2$ with $g_1 \in N_G(P)$
and $g_2 \in C_G(Z)$.
\end{lemma}

\begin{proof}
Observe that $P$~is a Sylow $p$-subgroup of $C_G(Z)$, and so all
Sylow $p$-subgroups of $C_G(Z)$ have Frattini subgroup~$Z$. Now set $R$ equal
to $P^g \cap C_G(Z)$. By assumption the Frattini subgroup of~$R$ is that
of~$P^g$, namely~$Z^g$. But $R$ is contained in a Sylow $p$-subgroup
of~$C_G(Z)$.  We conclude that $Z^g = Z$.

Therefore $P^g$ is itself a Sylow $p$-subgroup
of $C_G(Z)$, and so $P^g = P^h$ for some $h \in C_G(Z)$.
Take $g_1 = g h^{-1}$ and $g_2 = h$.
\end{proof}

\begin{lemma}
\label{lemma:Y}
Suppose that $P$ has odd exponent~$p^2$.
Then $Y^g = Y$ for every $g \in C_G(Z)$ such that $Y \leq P \cap P^g$.
Moreover, such $g$ factorize as $g_1 g_2$, where
$g_1$~lies in $C_G(Z) \cap N_G(P)$
and $g_2 \in C_G(Z) \cap N_G(Y)$ acts trivially on $Y/Z$.
\end{lemma}

\begin{proof}
Set $D_1$ equal to $C_G(Y)$, which contains $\Omega_1(P)$.  Since the centre
of~$P$ is cyclic and $\Omega_1(P)$ is maximal in~$P$, it follows that
$\Omega_1(P)$ is a Sylow $p$-subgroup of~$D_1$.  Now let $R$ be $P^g \cap D_1$,
the centralizer of $Y$~in $P^g$. Since $\Omega_1(P)$ has exponent~$p$,
so does~$R$. As $g$~centralizes $Z$, we deduce that $R$ is a maximal subgroup
of~$P^g$. But the only exponent~$p$ maximal subgroup of~$P^g$ is
$\Omega_1(P^g)$, which has centre $Y^g$.  So $g$~normalizes $Y$.

Now set $D_2 = \{h \in C_G(Z) \cap N_G(Y) \mid \text{$g$~acts trivially
on $Y/Z$}\}$\@. Then $P$~is a Sylow $p$-subgroup of~$D_2$, and $P^g$~is
too since $g$~lies in $C_G(Z) \cap N_G(Y)$. So $P^g = P^h$ for some
$h \in D_2$.  Take $g_1 = g h^{-1}$ and $g_2 = h$.
\end{proof}

\section{Stability conditions}
The following elementary reformulation of the usual stability
condition is not new, but does not appear to be widely known.

\begin{lemma}
\label{lemma:stableAfter}
Let $P$~be a Sylow $p$-subgroup of a finite group~$G$.  The cohomology class
$x \in H^*(P, \f[p])$ lies in $\Bild \Res^G_P$ if and only if
$x$~is an $N_G(P)$-invariant and
\[
\Cor^P_{P \cap P^g} \Res^{P^g}_{P \cap P^g} g^* (x) = 0 \quad
\text{for all $g \in G - N_G(P)$.}
\]
\end{lemma}

\begin{proof}
If $x$~comes from $H^*(G)$, then it certainly satisfies both conditions.
Conversely, observe that the conditions combined with the Mackey formula mean
that $\Res^G_P \Cor^G_P (x) = |N_G(P):P| x$.
\end{proof}

\noindent
Let~$P$ be a finite $p$-group.  The ring of universally stable elements~$I(P)$
was defined in~\cite{EvensPriddy} as the subring of $H^*(P,\f[p])$ given by
\[
I(P) = \bigcap_G \Bild \Res^G_P  \, ,
\]
where $G$ ranges over all finite groups with Sylow $p$-subgroup~$P$.
The following observation appears in~\cite{Ogawa:extraspecial}.
Recall that $O^p(G)$ is the subgroup generated by all $p'$-elements of~$G$.

\begin{lemma}
\label{lemma:Ogawa}
$I(P) \subseteq H^*(P)^{O^p(\Out(P))}$.
\end{lemma}

\begin{proof}
Pick any outer automorphism of order prime to~$p$, and lift it to an
automorphism $\alpha$~of the same order.  Let~$G$ be the semidirect
product $P \rtimes \langle \alpha \rangle$.
\end{proof}

\noindent
It is immediate that if $P$~is a Swan group then equality holds in
Lemma~\ref{lemma:Ogawa}.  The main result of this paper is:

\begin{theorem}
\label{theorem:main}
All extraspecial $p$-groups~$P$ apart from $D_8$~and $E$
are Swan groups.
If $P$ is $D_8$~or $E$, then the universally stable elements
for~$P$ are strictly contained in $H^*(P)^{O^p(\Out(P))}$.
\end{theorem}


\begin{proof}
See~\cite{MartinoPriddy:Swan} for a proof that $M(p^3)$ is a Swan group.
The groups $D_8$~and $E$ are treated in Lemma~\ref{lemma:exceptions} below.
So we may assume that $P$~satisfies the hypotheses of
Lemma~\ref{lemma:centralizerFrattini} and hence those of
Lemma~\ref{lemma:Z}.

Let $G$~be a finite group with extraspecial Sylow $p$-subgroup~$P$, and let
$x \in H^*(P)$ be an $N_G(P)$-invariant.  We show that the conditions of
Lemma~\ref{lemma:stableAfter} are satisfied.
Let $g$~be an element of~$G - N_G(P)$.
If $P^g \cap P$ is elementary abelian but not maximal in~$P$, then
corestriction from $P^g \cap P$ to~$P$ is zero: for corestriction from any
group $H$~to $H \times C_p$ is zero.

We may therefore assume that $P^g$, $P$ both contain $Z = \Phi(P)$.
Write $H$~for $P \cap P^g$ and $L$~for ${{}^g P} \cap P$.
By Lemma~\ref{lemma:Z} we deduce that $g$~normalizes $Z$. Moreover,
since $x$ is invariant under $N_G(P)$, we may in fact assume that
$g$~stabilizes $Z$.
If $P$~has odd exponent~$p^2$ we deduce further by Lemma~\ref{lemma:Y}
either that $g$~normalizes $Y$ and can be taken to act trivially on $Y/Z$,
or that $H \cap Y = L \cap Y = Z$.

So we can now apply Proposition~\ref{prop:Witt} and deduce that conjugation
$c_g \colon H \rightarrow L$ extends to an automorphism $\phi$~of $P$.
Then $\Res^{P^g}_H g^* = \Res^P_H \phi^*$, which means that
$\Cor^P_H \Res^{P^g}_H g^* = 0$ as $H$~is a proper subgroup of~$P$.
\end{proof}

\begin{lemma}
\label{lemma:exceptions}
If $P$~is $D_8$~or $E$, then the inclusion in Lemma~\ref{lemma:Ogawa} is
strict.
\end{lemma}

\begin{proof}
Let $G$~be $\GL_3(\f)$.  Then the upper triangular matrices with ones on
the diagonal form a Sylow $p$-subgroup isomorphic to~$P$.  We may take
$B_1 = \left(
\begin{smallmatrix} 1 & 0 & 0 \\ 0 & 1 & 1 \\ 0 & 0 & 1 \end{smallmatrix}
\right)$ and $C = \left(
\begin{smallmatrix} 1 & 0 & 1 \\ 0 & 1 & 0 \\ 0 & 0 & 1 \end{smallmatrix}
\right)$.  Pick $g = \left(
\begin{smallmatrix} 0 & 1 & 0 \\ 1 & 0 & 0 \\ 0 & 0 & 1 \end{smallmatrix}
\right)$.
Then $F = P \cap P^g$ equals $\langle B_1, C \rangle$, a maximal elementary
abelian subgroup.
Let $\beta, \gamma$ be the dual basis for~$F^*$.

Taking the first Chern class is an $\f$-linear monomorphism from $F^*$
to $H^2(F)$. Let~$\rho$ be an ordinary representation of~$P$ with
character~$\chi$, the sum of all $p^2$ linear characters.  These linear
characters restrict to~$F$ as scalar multiples of~$\beta$, each scalar
multiple being the image of $p$ characters.  So by the Whitney sum formula,
the total Chern class $c(\rho)$ restricts to~$F$ as follows:
\[
\Res^P_F c(\rho) = \prod_{\mu \in \f} (1 + \mu \beta)^p \, .
\]
This equals $1 - \beta^{p(p-1)}$. Set $\eta$ equal to $c_{p(p-1)}(\rho)$
in $H^{2p(p-1)}(P)$.  Then~$\eta$ lies in $H^*(P)^{O^p(\Out(P))}$,
since~$\chi$ is an invariant of $\Aut(P)$.
But $\Res^P_F(\eta) = -\beta^{p(p-1)}$ is distinct from
$g^* \Res^P_F(\eta) = - \gamma^{p(p-1)}$, so $\eta$~is not stable.
\end{proof}

\noindent
Suppose that $P$ is a $p$-group whose classifying space does not have
any transfer summands in its stable splitting. Then $P$~is a Swan group
by Theorem~3.5 of~\cite{MartinoPriddy:Swan}\@.  Martino and Priddy give one
counterexample to the converse: $M(p^3)$.

\begin{corollary}
\label{coroll:transfer}
Let $P$ be $2^{1+2n}_+$ with $n \geq 2$; or $p^{1+2n}_+$ with $p$~odd and
$n \geq 2$; or $p^{1+2n}_-$ with $p$~odd and $n \geq 1$. Then although~$P$
is a Swan group, the stable splitting of $BP$ involves a transfer summand.
\end{corollary}

\begin{proof}
$F = \langle B_1, \ldots, B_n, C \rangle$ is a self-centralising, maximal
elementary abelian subgroup. By Theorem~0.1 of~\cite{Priddy:Ln}, the
Steinberg summand $L(n+1)$ of~$BF$ is a transfer summand of~$BP$.
\end{proof}

\begin{remark}
\label{remark:Ogawa}
Ogawa's proof of Theorem~\ref{theorem:main} hinges around the following
claim: if $M$ is a maximal subgroup of an extraspecial $2$-group $P \not = D_8$,
then $P$~acts trivially on $H^*(M)$.  As observed in~\cite{Minh:proper}, the proof
of this claim in~\cite{Ogawa:extraspecial} uses inflation incorrectly.
We shall now see that the claim is false for~$2^{1+2n}_+$.

Let~$P$ be $D_8 * \cdots * D_8$, or $p^{1+2n}_+$ if $p$~is odd.
Let $M$~be $C_P(B_1)$, an index~$p$ subgroup of~$P$.  Then
$F = \langle B_1, \ldots, B_n, C \rangle$ is maximal elementary abelian
in $M$~and in $P$. Let $\beta_1, \ldots, \beta_n, \gamma$
be the dual basis of $F^*$, which we embed in $H^2(F)$ by taking the
first Chern class. Then $S(F^*) \subseteq H^*(F)$.
Define $\zeta \in H^{2p^{n-1}}(M)$ to be the Evens norm $N^M_F(\gamma)$,
and set $g = A_1$.  Let $U$~be the subspace of~$F^*$ spanned by
$\beta_2,\ldots,\beta_n$.  By standard properties of the Evens norm map
(see \cite[Ch.~6]{Evens:book}), we have:
\[
\Res^M_F(\eta) = \prod_{u \in U} (\gamma + u) \quad \text{and so} \quad
g^* \Res^M_F(\eta) = \prod_{u \in U} (\gamma + \beta_1 + u)
\, .
\]
These restrict to $\langle B_1, C \rangle$ as $\gamma^{p^{n-1}}$ and
$(\gamma + \beta_1)^{p^{n-1}}$ respectively.  So
$\Res^M_F(g^* \eta)$ differs from $\Res^M_F(\eta)$,
which means that $\eta$ is not invariant under~$A_1$.
\end{remark}

\bibliographystyle{abbrv}
\bibliography{ogawa}

\begin{thebibliography}{1}

\bibitem{Dieudonne:3rdEd}
J.~A. Dieudonn{\'e}.
\newblock {\em La g\'eom\'etrie des groupes classiques}.
\newblock Ergebnisse der Mathematik und ihrer Grenzgebiete, Band 5.
  Springer-Verlag, Berlin, third edition, 1971.

\bibitem{Evens:book}
L.~Evens.
\newblock {\em The cohomology of groups}.
\newblock Oxford Univ.\@ Press, Oxford, 1991.

\bibitem{EvensPriddy}
L.~Evens and S.~Priddy.
\newblock The ring of universally stable elements.
\newblock {\em Quart. J. Math. Oxford Ser. (2)}, 40(160):399--407, 1989.

\bibitem{Jacobson:BasicI}
N.~Jacobson.
\newblock {\em Basic algebra. {I}}.
\newblock W. H. Freeman and Company, New York, second edition, 1985.

\bibitem{MartinoPriddy:Swan}
J.~Martino and S.~Priddy.
\newblock On the cohomology and homotopy of {S}wan groups.
\newblock {\em Math. Z.}, 225(2):277--288, 1997.

\bibitem{Minh:proper}
P.~A. Minh.
\newblock Proper singularities of extraspecial $p$-groups for the mod-$p$
  cohomology functor.
\newblock {\em Comm. Algebra}, 25(3):965--971, 1997.

\bibitem{Ogawa:extraspecial}
Y.~Ogawa.
\newblock On the subring of universally stable elements in a mod-$2$ cohomology
  ring.
\newblock {\em Tokyo J. Math.}, 15(1):91--97, 1992.

\bibitem{Priddy:Ln}
S.~B. Priddy.
\newblock On characterizing summands in the classifying space of a group. {I}.
\newblock {\em Amer. J. Math.}, 112(5):737--748, 1990.

\end{thebibliography}

\end{document}